\documentclass[10pt]{amsart}
\usepackage{bm}
\usepackage{pstricks}

\newcommand{\wX}{\widetilde X}
\newcommand{\bA}{\mathbb A}
\newcommand{\bP}{\mathbb P}

\newcommand{\bQ}{\mathbb Q}
\newcommand{\bR}{\mathbb R}
\newcommand{\bC}{\mathbb C}

\newcommand{\dra}{\dashrightarrow}

\newcommand{\TOP}{\rm top}
\newcommand{\bottom}{\rm bot}

\DeclareMathOperator{\preimage}{\operatorname{preimage}}
\DeclareMathOperator{\Pic}{\operatorname{Pic}}
\DeclareMathOperator{\rank}{\operatorname{rank}}
\DeclareMathOperator{\image}{\operatorname{image}}

\def\OO{\mathcal O}

\theoremstyle{plain}
\newtheorem{theorem}{Theorem}[section]

\newtheorem{lemma}[theorem]{Lemma}
\newtheorem{corollary}[theorem]{Corollary}

\newtheorem{proposition}[theorem]{Proposition}

\theoremstyle{definition}

\newtheorem{definition}[theorem]{Definition}
\newtheorem{remark}[theorem]{Remark}

\newtheorem{example}[theorem]{Example}

\newtheorem*{acknowledgements}{Acknowledgments}

\author{Joe Rusinko}
\address{Department of Mathematics, University of Georgia, Athens GA
  30602, USA   email: rusinko@math.uga.edu}

\title[Equivalence of Mirror Families]{Equivalence of Mirror Families Constructed from Toric Degenerations of Flag Varieties}
\date\today

\begin{document}
\begin{abstract}Batyrev (et. al.) constructed a family of Calabi-Yau varieties using small toric degenerations of the full flag variety 
$G/B$. They conjecture this family to be mirror to generic anti-canonical hypersurfaces  in $G/B$.  Recently Alexeev 
and Brion, as a part of their work on toric degenerations of spherical varieties, have constructed many degenerations 
of $G/B$.   For any such degeneration we construct a family of varieties, which we prove coincides with Batyrev's in 
the small case.  We prove that any two such families are birational, thus proving that mirror families are independent 
of the choice of degeneration.  The birational maps involved are closely related to Berenstein and Zelevinsky's geometric 
lifting of tropical maps to maps between totally positive varieties.

\end{abstract}
\maketitle

\section{Introduction}
\label{intro}
Generic elements of the anti-canonical class of a flag variety $G/P$ are smooth Calabi-Yau varieties.  We restrict our
attention to $G$ of type $A_n$ (i.e. $G=SL_{n+1}$), and work over the ground field $\bC$.  For any parabolic subgroup $P$,
Batyrev, Ciocan-Fontanine, Kim, and van Straten \protect\cite{partialflag} use Gonciulea and Lakshmibai's 
\protect\cite{Lakshmibai} toric degeneration of $G/P$, to construct a family of CYs which they conjecture to be mirror to the 
anti-canonical hypersurfaces in $G/P$.  In particular, they degenerate the pair ($G/B$,$-K_{G/B}$) to a toric pair
($X_{\Delta}$,$-K_{X_\Delta}$) corresponding to a certain polytope $\Delta$.  The polytope $\Delta$ has a dual 
polytope $\Delta^*$ corresponding to a toric pair ($X_{\Delta^*}$,$-K_{X_{\Delta^*}}$).  They
define a subfamily of $|-K_{X_{\Delta^*}}|$ whose coefficients satisfy ``box equations''. They conjecture that 
minimal crepant resolutions of generic elements of this subfamily are mirror to generic anti-canonical hypersurfaces 
in $G/B$. As evidence, they prove that the varieties have the same associated hypergeometric series 
\protect\cite[5.2.1]{partialflag}.

Following Caldero's work in \protect\cite{cald}, Alexeev and Brion have constructed many different degenerations
of spherical varieties \cite{ab}.  These constructions include degenerations of the flag variety indexed by reduced 
decompositions $\overline{\omega_0}$ of $\omega_0$, the longest element in the Weyl group of $G$.  In particular ($G/B$,$-K_{G/B}$) 
degenerates to a toric pair ($X_{\Delta},\OO_{X_{\Delta}}(1)$) 
corresponding to a rational polytope $\Delta = \Delta(\overline{\omega_0})$ \protect\cite{ab}.  We call these degenerations
\emph{string degenerations}.  The fibers in the limit vary greatly depending on the choice of reduced decomposition.  For example,
the polytope's F-vector (number of faces 
of different dimensions), the rank of the Picard group of $X_{\Delta}$, and the dimension of the linear system
$|\OO_{X_{\Delta^*}}(1)|$ may all differ. 

We would like to construct a mirror family using any of the string degenerations and see how they
differ depending on the choice of degeneration.  For a small  toric degeneration (see Definition \ref{Def:small}) 
Batyrev proposed such a family \protect\cite{B_mirror}. 

We prove that $\Delta^*(\overline{\omega_0})$ is integral, so it corresponds to a 
toric variety $X_{\Delta^*}$ with an ample line bundle $\OO_{X_{\Delta^*}}(1)$.  For any choice of reduced composition 
$\overline{\omega_0}$, we define $F_{\overline{\omega_0}}$ as a subfamily of  
$|\OO_{X_{\Delta^*}}(1)|$ whose coefficients satisfy a set of relations we call \emph{combinatorial box equations}. 
We prove that this family coincides with the family in Batyrev's construction when the string degeneration is small.

Using Berenstein and Zelevinsky's concept of geometric lifting \protect\cite{bz}, we construct birational maps between 
any two families $F_{\overline{\omega_0}}$ and $F_{\overline{\omega_0}'}$.    By the work of Batyrev and Kontsevich
\protect\cite{bbetti, Kont}, this shows that smooth Calabi-Yau  varieties birational to the generic elements of
$F_{\overline{\omega_0}}$ all have the same Hodge numbers. 

The organization of the paper is as follows:

In section \protect\ref{sec:preliminaries} we  fix notation and describe
the string degenerations of $G/B$. 
We review Gleizer and Postnikov's description of the string cone \protect\cite{gp} and use it to define a
family of varieties $F_{\overline{\omega_0}}$.
  
Every pair of reduced decompositions of $\omega_0$ is connected by a sequence of braid moves.  In section
\protect\ref{sec:inequalities}, we classify how the inequalities defining the polytopes $\Delta$ change under a braid
move.  This allows us to prove that the polytopes $\Delta^*$ are integral, and will
help us in proving that $F_{\overline{\omega_0}}$ and $F_{\overline{\omega_0}'}$ are 
birational.

In section \protect\ref{sec:birational}, we define a birational map between the families $F_{\overline{\omega_0}}$ 
and $F_{\overline{\omega_0}'}$.  This map is closely related to Berenstein and Zelevinsky's geometric 
lifting of piecewise linear maps to maps between totally positive varieties \protect\cite{bz}.   

In section \protect\ref{sec:small}, we prove that the families
$F_{\overline{\omega_{0}}}$ coincide with the families that Batyrev constructed when
the string degeneration is small. Finally, we give an example of a string degeneration which is not small.
\begin{acknowledgements}
The author would like to thank his advisor Valery Alexeev for posing this interesting problem, and for his helpful
guidance and advice.  The author would also like to thank Robert Varley
for many helpful conversations. 
\end{acknowledgements}

\section{Preliminaries}
\label{sec:preliminaries}
We begin by introducing notation that will be used throughout this paper.  
Working over $\bC$, let $G$ be the algebraic group $SL_{n+1}$, $B$ the Borel subgroup of upper triangular
matrices, $U$ the subgroup of upper uni-triangular matrices, and $T$ the maximal torus of
diagonal matrices. One has $B=TU$.  Let $\Phi$ be the root system of ($G$,$T$), with  $\Phi^+$  the subset of positive roots.
We denote $\alpha_1,...,\alpha_n$ the corresponding simple roots.  

Let $W$ be the Weyl group of ($G$,$T$) with $s_1,...,s_n$  the reflections associated to 
$\alpha_1,\dots,\alpha_n$.  Note that $W \cong S_{n+1}$ so we can view $s_i$ as the transposition  exchanging
$i$ and $i+1$.  These transpositions define a length function $l$ on $W$. Let $\omega_0$ be the unique element
of maximal length corresponding to the permutation sending $1$ to $n+1$, $2$ to $n$ and so forth.  There are
$N=l(\omega_0)=\frac{n(n+1)}{2}$ positive roots, which is the dimension of the flag variety $G/B$.

We say that a decomposition $\overline{\omega} = s_{i_1}s_{i_2}\dots s_{i_m}$ of $\omega$ is 
reduced if $m=l(\omega)$.   For every reduced decomposition of $\omega_0$,
the string degenreation is a flat degeneration of ($G/B$,$-K_{G/B}$) to a toric variety $X_{\Delta}$ 
and a corresponding $\bQ$-line bundle $\OO_{X_{\Delta}}(1)$ \cite{ab,cald}.  By standard methods 
(see \protect\cite{fulton,toricmirror}) one can associate to ($X_{\Delta}$,$\OO_{X_{\Delta}}(1)$) a polytope 
$\Delta=\Delta(\overline{\omega_0})$ and its dual polytope $\Delta^*=\Delta^*(\overline{\omega_0})$.  For any
reduced decomposition, $\Delta$ is a rational polytope in $\bR^N$, thus $X_{\Delta}$ is 
a toric variety of $\dim N$.
For $n \le 4$ and any choice of reduced 
decomposition, Alexeev and Brion  verified that $\Delta$ is integral  \protect\cite[Ex. 5.6]{ab}.  This implies that the  
$\bQ$-line bundle is actually a line bundle.   In Section \protect\ref{sec:inequalities} we prove that for any $n$ and
any choice of reduced decomposition, $\Delta^*$ is 
an integral polytope and thus corresponds to a toric variety $X_{\Delta^*}$ with ample line bundle 
$\OO_{X_{\Delta^*}}(1)$.

The polytope $\Delta$ is the intersection of Berenstein-Zelevinsky's string cone and a polyhedral
cone called the $\lambda$-cone \protect\cite[Thm 1.1]{ab}. Explicit
inequalities describing the $\lambda$-cone can be found in \protect\cite[1]{littleman} and \protect\cite[Thm 1.1]{ab}. 
Inequalities for the string cone are given in \protect\cite[3.10]{bz} and \protect\cite[1]{littleman}. 
A combinatorial description for them is given in \protect\cite[Cor.5.8]{gp}.  

Gleizer and Postnikov's description of the string cone is crucial to what follows
so we review it here. (Note: in this review the order in which the $s_i$ appear is reversed from 
that of the original.)

Let $\overline{\omega_0}$=$s_{i_1}s_{i_2}\dots s_{i_N}$ be a reduced decomposition of $\omega_0$.
They represent this decomposition with a string diagram described below:
Start with $n+1$ strings at the top of the diagram.  Move down the diagram exchanging the string in the
$i_1$ column with the string in the $i_1$+1 column.  Continuing down the diagram exchange the string
in the $i_2$ column with the string in the $i_2 +1$ column.  Continue this process until you have made
$N$ exchanges.  The strings will have reversed their orders when they 
reach the bottom of the diagram.  

Label the intersection points $t_1, t_2,\dots,t_N$ from top to bottom.  Label the strings $U_1$ through $U_{n+1}$
on the top from left to right.  Mark vertices $u_1$ through $u_{n+1}$(resp. $b_1$ 
through $b_{n+1}$) as the upper (resp. lower) ends of the strings $U_1$ through $U_{n+1}$.  

\begin{example} Here is a picture of the string diagram for $n=3$ and reduced word decomposition 
$\overline{\omega_0}$ = $s_3s_1s_2s_1s_3s_2$: 

\begin{center}
\begin{pspicture}(0,-.3)(5,4.2)
\psset{unit=.2}
\psline(0,19)(0,15)(4,14)(8,13)(8,12)(12,11)(16,10)(16,6)(20,5)(24,4)(24,0)
\psline(8,19)(8,15)(4,14)(0,13)(0,9)(4,8)(8,7)(8,3)(12,2)(16,1)(16,0)
\psline(16,19)(16,18)(20,17)(24,16)(24,6)(20,5)(16,4)(16,3)(12,2)(8,1)(8,0)
\psline(24,19)(24,18)(20,17)(16,16)(16,12)(12,11)(8,10)(8,9)(4,8)(0,7)(0,0)
\uput[90](0,19){$u_1$}
\uput[90](8,19){$u_2$}
\uput[90](16,19){$u_3$}
\uput[90](24,19){$u_4$}
\uput[90](20,17){$t_1$}
\uput[90](4,14){$t_2$}
\uput[90](12,11){$t_3$}
\uput[90](4,8){$t_4$}
\uput[90](20,5){$t_5$}
\uput[90](12,2){$t_6$}
\uput[270](0,0){$b_4$}
\uput[270](8,0){$b_3$}
\uput[270](16,0){$b_2$}
\uput[270](24,0){$b_1$}
\end{pspicture}
\end{center}
\end{example}
\begin{definition} (\protect\cite{gp})
Pick a number $k$ from $1$ to $n$. Form an oriented graph on the string diagram as follows.   Orient strings $U_1$
through $U_k$ upward and the rest of the strings downward.  A \emph{rigorous path} is an oriented path from 
$b_k$ to $b_{k+1}$ not meeting any vertex more than once, and avoiding the two bad fragments pictured below.
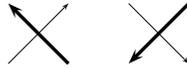
\begin{center}
\begin{pspicture}(0,0)(3,1)
\psset{unit=.4}
\psline[linewidth=.2 pt]{->}(0,0)(2,2)
\psline[linewidth=1.3 pt]{<-}(0,2)(2,0)
\psline[linewidth=1.3 pt]{<-}(4,0)(6,2)
\psline[linewidth=.2 pt]{->}(4,2)(6,0)
\end{pspicture}
\end{center}
Here the thin line denotes a string, and the thick line indicates the rigorous path.
\end{definition}
Each rigorous path defines an inequality as follows: 
\[
\Sigma a_ix_i \ge 0  \text{ where }a_i =
\begin{cases}
 1  &\text{if the path switches from }  U_i \text{ to } U_j \text{ and } i<j  \\
 -1 &\text{if the path switches from }  U_i \text{ to } U_j \text{ and } i>j  \\
\end{cases}
\]
We refer to these inequalities as string inequalities.
\begin{example}
For $k=1$ in the string diagram pictured above, the rigorous path $b_1 \mapsto t_5 \mapsto  t_3 \mapsto t_4 
\mapsto t_6 \mapsto b_2$ corresponds to the inequality: $x_3 \ge x_4$. 
\end{example}
Gleizer and Postnikov prove that the string cone is the collection of points in $\bR^N$ satisfying the string inequalities \protect\cite[Cor. 5.8]{gp}.
\begin{lemma}
The $\lambda$-cone is the collection of points in $\bR^N$ satisfying the following inequalities which
we call $\lambda$-inequalities: \\
\[
\lambda_i  :   x_i \leq 2 + \sum_{j > i}  c_j x_j
\]
for $i$ from $1$ to $N$, where
\[
c_j=
\begin{cases}
  1 &\text{if the vertex }  t_j  \text{ is one column to the
right or left of } t_i  \\
 -2 &\text{if the vertex for }  t_j \text{ is in the same column as } t_i  \\
 0 &\text{otherwise.} \\
\end{cases}
\]
\end{lemma}
This is simply a way of visualising the Lie-algebraic definition of the $\lambda$-cone  (see \protect\cite[Thm 1.1]{ab}), 
and of associating a $\lambda$-inequality to each intersection point in the string diagram.
\begin{example}
For the string diagram pictured above we have \[\lambda_2 : x_2 \le 2 + x_3 - 2x_4 + x_6.\]
\end{example}
\begin{definition}  Write the defining inequalities of $\Delta(\overline{\omega_0})$ as
$\sum_i m_{i_d}x_i \le $(0 or 2), and define $M_d =(m_{1_d},\cdots,m_{N_d})$. Let 
$T^{M_d}:=t_1^{m_{1_d}}t_2^{m_{2_d}}\cdots t_N^{m_{N_d}}$ and define \emph{$BF_{\overline{\omega_0}}$} 
as the family of hypersurfaces in the torus $T$=Spec$\bC[t_1,t_1^{-1},\dots,t_N,t_N^{-1}]$ given by the equations:
\[f_{\overline{\omega_0},a}:=  1 - \sum_{d} a_dT^{M_d}=0\] where $d$ runs over all defining inequalities of $\Delta$.
We will refer to individual hypersurfaces as $Z_{\overline{\omega_0},a}$, where $a=(a_1,\cdots, a_r)$ denotes the 
coefficient vector. 
\end{definition}
\begin{example}
The $\lambda_2$-inequality above would correspond to the term $a_{\lambda_2}t_2t_3^{-1}t_4^2t_6^{-1}$.
\end{example}
Note that $BF_{\overline{\omega_0}}$ is a family over the affine space $\bA^r$.
We will see
in Section \protect\ref{sec:inequalities} that it is a subfamily of the linear system 
$|\OO_{X_{\Delta^*}}(1)|$ intersected with the torus.

We restrict our attention to a subfamily  whose coefficients satisfy a set of 
relations we call \emph{combinatorial box equations}.  We define these relations in terms of the bounded regions of the 
string diagram which we call \emph{boxes}. Every box of the string diagram is obviously bounded above by a vertex 
$t_{\TOP}$ and directly below by a vertex $t_{\bottom}$.  Let $T^{\lambda_{\TOP}}$ (resp. $T^{\lambda_{\bottom}}$) 
be the monomial in $f_{\overline{\omega_0},a}$ corresponding to the $\lambda$-inequality associated to $t_{\TOP}$ (resp. 
$t_{\bottom}$).

\begin{definition}
For every two string inequalities $p_1$, $p_2$ with corresponding monomials $T^{p_1}$ and $T^{p_2}$ satisfying 
the following \emph{box conditions}:
\begin{enumerate}
\item there exists a box with corresponding monomials $T^{\lambda_{\TOP}}$ and 
$T^{\lambda_{\bottom}}$ such that $T^{p_1}$ $T^{\lambda_{\TOP}}$ = $T^{p_2}$ $T^{\lambda_{\bottom}}$,
\item the $t_{\TOP}$ degree of $T^{p_1} = -1$,
\item the $t_{\bottom}$ degree of $T^{p_2} = 1$,
\end{enumerate}
we define an equation $a_{p_1}a_{\lambda_{\TOP}} = a_{p_2}a_{\lambda_{\bottom}}$.  We call the collection
of all such equations the \emph{combinatorial box equations}.
\end{definition}
\begin{definition}
Let $P_{\overline{\omega_0}}$ be the subscheme defined by the combinatorial box equations in the torus Spec 
$\bC[a_1,a_1^{-1},\cdots,a_r,a_r^{-1}]$.
\end{definition}
\begin{definition}
Let  \emph{$F_{\overline{\omega_0}}$} be the subfamily of  $BF_{\overline{\omega_0}}$   whose coefficients are 
non-zero and satisfy the combinatorial box equations.  Note that $F_{\overline{\omega_0}}$ is a family over the
base $P_{\overline{\omega_0}}$.
\end{definition}
It is a standard fact that any two reduced decompositions of $\omega_0$  are connected by a finite sequence of the 
following two braid moves.
\begin{enumerate}
\item 2-move: exchanges ($s_i$, $s_j$) with ($s_j$, $s_i$) where $|i-j| > 1$
\item 3-move: exchanges ($s_i,s_j,s_i$) with ($s_j,s_i,s_j$) where $|i-j| = 1$
\end{enumerate}

A 3-move from $\overline{\omega_0}$ to $\overline{\omega_0}'$ fixes all of the string diagram except for 
one box.
\begin{definition}
Given a 3-move from $\overline{\omega_0}$ to $\overline{\omega_0}'$,  label the box in the string diagram for $\overline{\omega_0}$
bounded by $t_i, t_j$ and $t_k$, $R$ and its image $R'$. 
We say that $R$ is of type $R121$ if it is of the form ($s_i,s_{i+1},s_i$).  Otherwise we say it is of type $R212$. 
We refer to these regions as braid regions.  Note that if $R$
is of type $R121$ then $R'$ must be of type $R212$ and vice versa.
\end{definition}
The map between string cones differing by a single braid move, defined by Berenstein and Zelevinsky in \protect\cite[Thm 2.7]{bzan},
restricts to a map between polytopes $\Delta(\overline{\omega_0})$ and $\Delta(\overline{\omega_0}')$ .

They prove that the map fixes all of the coordinates 
except for the two or three corresponding to those being exchanged in the braid move.  For those coordinates they 
have the following maps:
\begin{enumerate}
\item 2-move: ($x_i,x_j$) $\rightarrow$ ($x_j,x_i$)
\item 3-move: ($x_i,x_j,x_k$) $\rightarrow$ (max($x_k$,$x_j - x_i$),$x_i + x_k$,min($x_i$,$x_j - x_k$)).
\end{enumerate}

\section{Classification of how $\Delta$ Changes Under 3-Move}
\label{sec:inequalities}
In this section we classify the defining inequalities of $\Delta$ and how those inequalities
change under a 3-move.  From these results we prove that the dual polytope $\Delta^*$ is integral.

Only three coordinates change under a 3-move.  In what follows we will refer to those coordinates as
$x_i$, $x_j$, and $x_k$ (note that $i, j,$ and $k$ are consecutive integers). We call the product of the $t_i$, $t_j$, 
and $t_k$ variables inside of $T^{M_d}$ a \emph{monomial piece}, and refer to this product as $\tau^{M_d}$.  
We will refer to the monomial pieces associated to 
string (resp. $\lambda$) inequalities as string (resp. $\lambda$) monomial pieces.

Since each inequality corresponds to a monomial, we can classify all possible string and $\lambda$-inequalities
in terms of their monomial pieces.
\begin{subsection}{Classification of $\lambda$-Monomial Pieces}
\begin{theorem}
The following is a complete classification of $\lambda$-monomial pieces before and after a 3-move.\\
$t_it_j^{-1}t_k^2 \leftrightarrow t_it_j^{-1}t_k^2 \qquad  t_jt_k^{-1} \leftrightarrow t_jt_k^{-1} \qquad  t_k  
\leftrightarrow t_k \qquad  t_i^{-1}t_k^{-1} \leftrightarrow t_j^{-1} \\ t_i^2t_j^{-1}t_k^2 \leftrightarrow 
t_i^{-1}t_j^2t_k^{-1}$ \\
An arrow indicates which monomial pieces are exchanges under the braid move.
\end{theorem}
\begin{proof}
As above, we fix $t_i,t_j$ and $t_k$ as the coordinates on which the 3-move occurs.  We classify
the monomial pieces for the inequalities $\lambda_l$ for all values of $l$.
For $l > k$ the monomial for $\lambda_l$ doesn't have any terms from the braid move section.  
This means that the monomial piece is $1$ and it doesn't change after a braid move.

For $\lambda_i$, $\lambda_j$, and $\lambda_k$ we have the following monomial pieces regardless if $R$ is of type $R121$ 
or $R212$: $\lambda_i  :  t_it_j^{-1}t_k^2$  \qquad $\lambda_j  :  t_jt_k^{-1}$   \qquad $\lambda_k  :  t_k   $ \\

If $l < i$, then the $\lambda_l$ monomial piece only depends on which column $t_l$ lies (as pictured below).
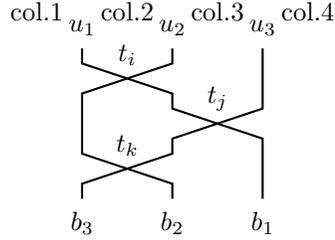
\begin{center}
\begin{pspicture}(-.7,-.3)(7,2.5)
\psset{unit=.2}
\psline(0,10)(0,9)(3,8)(6,7)(6,6)(9,5)(12,4)(12,0)
\psline(6,10)(6,9)(3,8)(0,7)(0,3)(3,2)(6,1)(6,0)
\psline(12,10)(12,6)(9,5)(6,4)(6,3)(3,2)(0,1)(0,0)
\uput[90](0,10){$u_1$}
\uput[90](6,10){$u_2$}
\uput[90](12,10){$u_3$}
\uput[90](3,8){$t_i$}
\uput[90](9,5){$t_j$}
\uput[90](3,2){$t_k$}
\uput[270](0,0){$b_3$}
\uput[270](6,0){$b_2$}
\uput[270](12,0){$b_1$}
\uput[90](-3,11){$\rm col. $1}
\uput[90](3,11){$\rm col. $2}
\uput[90](9,11){$\rm col. $3}
\uput[90](15,11){$\rm  col. $4}
\end{pspicture}
\end{center}
\begin{enumerate}
\item  If $t_l$ is in column $1$ we have the monomial piece  $t_i^{-1}t_k^{-1} $ before, and $t_j^{-1}$ after the 
braid move.
\item  If $t_l$ is in column $2$ we have the monomial piece $t_i^2t_j^{-1}t_k^2 $ before, and $ t_i^{-1}t_j^2t_k^{-1}$ 
after the braid move.
\item  If $t_l$ is in column $3$ we have the monomial piece $t_i^{-1}t_j^2t_k^{-1} $ before, and $t_i^2t_j^{-1}t_k^{2}$ 
after the braid move.
\item  If $t_l$ is in column $4$ we have the monomial piece $t_j^{-1} $ before, and $ t_i^{-1}t_k^{-1}$ 
after the braid move.
\end{enumerate}
Note that for $R$ of type $R212$ we get the same collection of monomial pieces.  
\end{proof}
\end{subsection}
\begin{subsection}{Classification of String Monomial Pieces}
\begin{theorem}
The following is a complete classification of the string monomial pieces, and how they are exchanged under a type 
$B_2$ braid move. The $\&$ sign indicates that a pair of monomial pieces must occur together.
\begin{align*}
t_k^{-1} \leftrightarrow &t_i^{-1} \text{ }\&\ \text{ } t_j^{-1}t_k \\
t_i \leftrightarrow &t_k  \text{ }\&\ \text{ } t_i^{-1}t_j \\
t_it_k \leftrightarrow &t_j \\
t_i^{-1}t_k^{-1} \leftrightarrow &t_j^{-1} \\
t_it_j^{-1}t_k \leftrightarrow &t_i^{-1}t_jt_k^{-1} \\
\end{align*}
\end{theorem}
\begin{proof}
We would like to classify all of the string inequalities in terms of their monomial pieces.  By symmetry
we can assume that $R$ is of type $R121$.  If a string monomial piece is anything other than $1$, the corresponding rigorous 
path $p$ must change strings in the region $R$. Any oriented path entering $R$ more than once must intersect itself so 
it cannot be rigorous. Therefore, we may assume that a rigorous path only enters and exits $R$ once. There are $6*5=30$ 
possible pairs of entry and exit points for $p$.  
If $p$ enters $R$ at the vertex $b_e$, then all strings $U_f$ for $f \ge e$ must be oriented upward.  Therefore $p$
cannot exit $R$ through $b_f$ where $f < e$.  Similarly if $p$ enters at $u_e$ it may not exit through $u_f$ for $f > e$. 
This eliminates $6$ pairs of entry and exit points.  By the same reasoning no rigorous path which changes
strings on $R$ can enter at $b_3$ (resp. $u_1$) and exit at $u_3$ (resp. $b_1$). What remains are $22$ entry and exit pairs.   

The following three diagrams correspond to cases of entry and exit points in which there is a unique rigorous path before
and after the braid move.
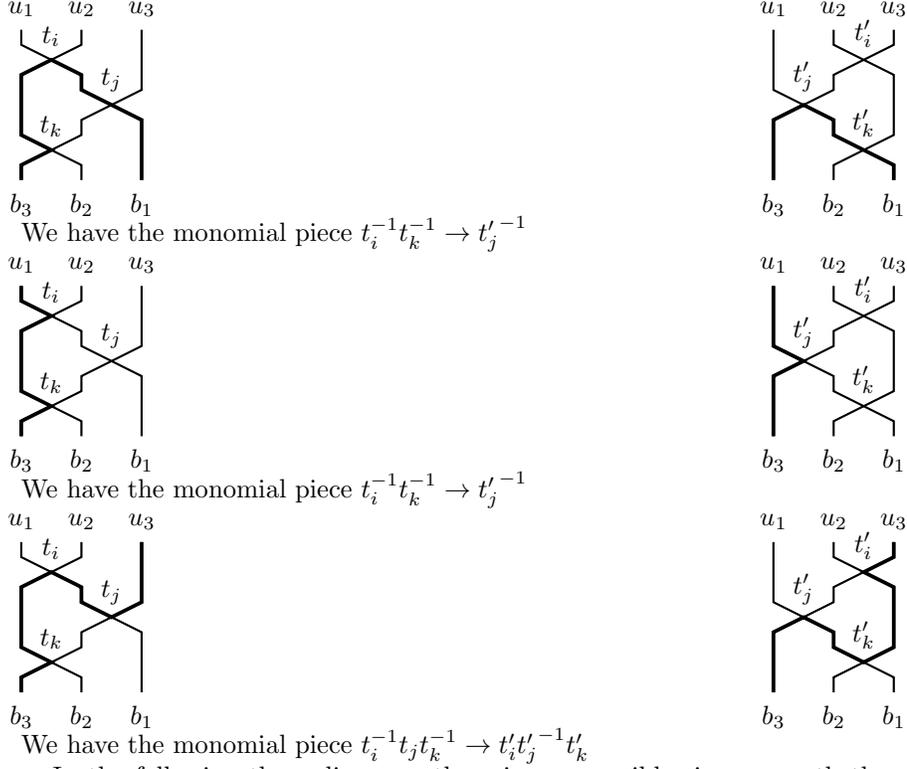
\begin{center}
\begin{pspicture}(0,-.15)(20,1.9)
\psset{unit=.2}
\psline(0,10)(0,9)(2,8)(4,7)(4,6)(6,5)(8,4)(8,0)
\psline(4,10)(4,9)(2,8)(0,7)(0,3)(2,2)(4,1)(4,0)
\psline(8,10)(8,6)(6,5)(4,4)(4,3)(2,2)(0,1)(0,0)
\uput[90](0,10){$u_1$}
\uput[90](4,10){$u_2$}
\uput[90](8,10){$u_3$}
\uput[90](2,8){$t_i$}
\uput[90](6,5){$t_j$}
\uput[90](2,2){$t_k$}
\uput[270](0,0){$b_3$}
\uput[270](4,0){$b_2$}
\uput[270](8,0){$b_1$}
\psline[linewidth=.25](8,0)(8,4)(6,5)(4,6)(4,7)(2,8)(0,7)(0,3)(2,2)(0,1)(0,0)

%2 212 piece
\psline(50,10)(50,6)(52,5)(54,4)(54,3)(56,2)(58,1)(58,0)
\psline(54,10)(54,9)(56,8)(58,7)(58,3)(56,2)(54,1)(54,0)
\psline(58,10)(58,9)(56,8)(54,7)(54,6)(52,5)(50,4)(50,0)
\uput[90](50,10){$u_1$}
\uput[90](54,10){$u_2$}
\uput[90](58,10){$u_3$}
\uput[90](56,8){$t_i'$}
\uput[90](52,5){$t_j'$}
\uput[90](56,2){$t_k'$}
\uput[270](50,0){$b_3$}
\uput[270](54,0){$b_2$}
\uput[270](58,0){$b_1$}
\psline[linewidth=.25](58,0)(58,1)(56,2)(54,3)(54,4)(52,5)(50,4)(50,0)
\end{pspicture}
\end{center}
We have the monomial piece $t_i^{-1}t_k^{-1} \rightarrow {t'_j}^{-1}$

\begin{center}
\begin{pspicture}(0,-.15)(20,1.9)
\psset{unit=.2}
\psline(0,10)(0,9)(2,8)(4,7)(4,6)(6,5)(8,4)(8,0)
\psline(4,10)(4,9)(2,8)(0,7)(0,3)(2,2)(4,1)(4,0)
\psline(8,10)(8,6)(6,5)(4,4)(4,3)(2,2)(0,1)(0,0)
\uput[90](0,10){$u_1$}
\uput[90](4,10){$u_2$}
\uput[90](8,10){$u_3$}
\uput[90](2,8){$t_i$}
\uput[90](6,5){$t_j$}
\uput[90](2,2){$t_k$}
\uput[270](0,0){$b_3$}
\uput[270](4,0){$b_2$}
\uput[270](8,0){$b_1$}
\psline[linewidth=.25](0,0)(0,1)(2,2)(0,3)(0,7)(2,8)(0,9)(0,10)

%2 212 piece
\psline(50,10)(50,6)(52,5)(54,4)(54,3)(56,2)(58,1)(58,0)
\psline(54,10)(54,9)(56,8)(58,7)(58,3)(56,2)(54,1)(54,0)
\psline(58,10)(58,9)(56,8)(54,7)(54,6)(52,5)(50,4)(50,0)
\uput[90](50,10){$u_1$}
\uput[90](54,10){$u_2$}
\uput[90](58,10){$u_3$}
\uput[90](56,8){$t_i'$}
\uput[90](52,5){$t_j'$}
\uput[90](56,2){$t_k'$}
\uput[270](50,0){$b_3$}
\uput[270](54,0){$b_2$}
\uput[270](58,0){$b_1$}
\psline[linewidth=.25](50,0)(50,4)(52,5)(50,6)(50,10)
\end{pspicture}
\end{center}
We have the monomial piece $t_i^{-1}t_k^{-1} \rightarrow {t'_j}^{-1}$

\begin{center}
\begin{pspicture}(0,-.15)(20,1.9)
\psset{unit=.2}
\psline(0,10)(0,9)(2,8)(4,7)(4,6)(6,5)(8,4)(8,0)
\psline(4,10)(4,9)(2,8)(0,7)(0,3)(2,2)(4,1)(4,0)
\psline(8,10)(8,6)(6,5)(4,4)(4,3)(2,2)(0,1)(0,0)
\uput[90](0,10){$u_1$}
\uput[90](4,10){$u_2$}
\uput[90](8,10){$u_3$}
\uput[90](2,8){$t_i$}
\uput[90](6,5){$t_j$}
\uput[90](2,2){$t_k$}
\uput[270](0,0){$b_3$}
\uput[270](4,0){$b_2$}
\uput[270](8,0){$b_1$}
\psline[linewidth=.25](0,0)(0,1)(2,2)(0,3)(0,7)(2,8)(4,7)(4,6)(6,5)(8,6)(8,10)

%2 212 piece
\psline(50,10)(50,6)(52,5)(54,4)(54,3)(56,2)(58,1)(58,0)
\psline(54,10)(54,9)(56,8)(58,7)(58,3)(56,2)(54,1)(54,0)
\psline(58,10)(58,9)(56,8)(54,7)(54,6)(52,5)(50,4)(50,0)
\uput[90](50,10){$u_1$}
\uput[90](54,10){$u_2$}
\uput[90](58,10){$u_3$}
\uput[90](56,8){$t_i'$}
\uput[90](52,5){$t_j'$}
\uput[90](56,2){$t_k'$}
\uput[270](50,0){$b_3$}
\uput[270](54,0){$b_2$}
\uput[270](58,0){$b_1$}
\psline[linewidth=.25](50,0)(50,4)(52,5)(54,4)(54,3)(56,2)(58,3)(58,7)(56,8)(58,9)(58,10)
\end{pspicture}
\end{center}
We have the monomial piece $t_i^{-1}t_jt_k^{-1} \rightarrow t'_i{t'_j}^{-1}t'_k$

In the following three diagrams there is one possible rigorous path through $R$ before the braid move.  After the braid
move there are two possible rigorous paths through $R'$.  Since the string cone is the collectin of points satisfying all possible string
inequalities there must be two inequalities (and hence monomial pieces) in the image.  In what follows we group these
monomial pieces together.

% page 1 # 1
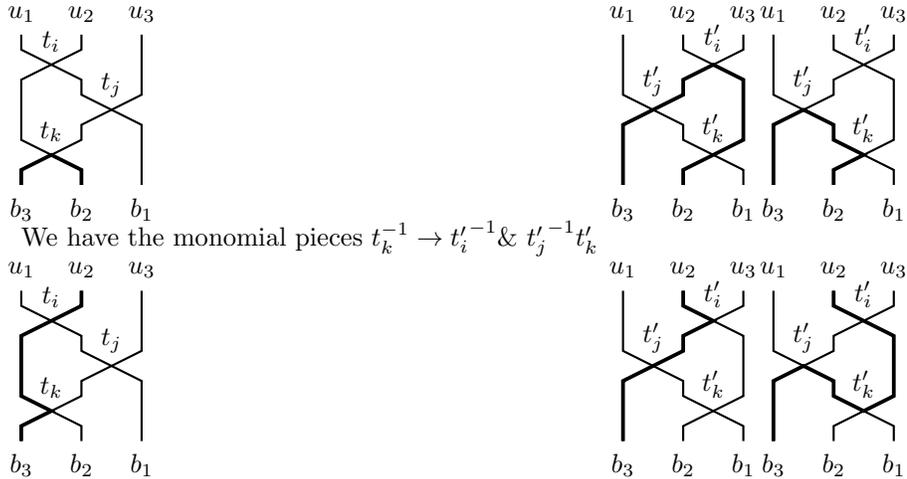
\begin{center}
\begin{pspicture}(0,-.15)(20,1.9)
\psset{unit=.2}
\psline(0,10)(0,9)(2,8)(4,7)(4,6)(6,5)(8,4)(8,0)
\psline(4,10)(4,9)(2,8)(0,7)(0,3)(2,2)(4,1)(4,0)
\psline(8,10)(8,6)(6,5)(4,4)(4,3)(2,2)(0,1)(0,0)
\uput[90](0,10){$u_1$}
\uput[90](4,10){$u_2$}
\uput[90](8,10){$u_3$}
\uput[90](2,8){$t_i$}
\uput[90](6,5){$t_j$}
\uput[90](2,2){$t_k$}
\uput[270](0,0){$b_3$}
\uput[270](4,0){$b_2$}
\uput[270](8,0){$b_1$}
\psline[linewidth=.25](4,0)(4,1)(2,2)(0,1)(0,0)

%1 212 piece
\psline(40,10)(40,6)(42,5)(44,4)(44,3)(46,2)(48,1)(48,0)
\psline(44,10)(44,9)(46,8)(48,7)(48,3)(46,2)(44,1)(44,0)
\psline(48,10)(48,9)(46,8)(44,7)(44,6)(42,5)(40,4)(40,0)
\uput[90](40,10){$u_1$}
\uput[90](44,10){$u_2$}
\uput[90](48,10){$u_3$}
\uput[90](46,8){$t_i'$}
\uput[90](42,5){$t_j'$}
\uput[90](46,2){$t_k'$}
\uput[270](40,0){$b_3$}
\uput[270](44,0){$b_2$}
\uput[270](48,0){$b_1$}
\psline[linewidth=.25](44,0)(44,1)(46,2)(48,3)(48,7)(46,8)(44,7)(44,6)(42,5)(40,4)(40,0)

%2 212 piece
\psline(50,10)(50,6)(52,5)(54,4)(54,3)(56,2)(58,1)(58,0)
\psline(54,10)(54,9)(56,8)(58,7)(58,3)(56,2)(54,1)(54,0)
\psline(58,10)(58,9)(56,8)(54,7)(54,6)(52,5)(50,4)(50,0)
\uput[90](50,10){$u_1$}
\uput[90](54,10){$u_2$}
\uput[90](58,10){$u_3$}
\uput[90](56,8){$t_i'$}
\uput[90](52,5){$t_j'$}
\uput[90](56,2){$t_k'$}
\uput[270](50,0){$b_3$}
\uput[270](54,0){$b_2$}
\uput[270](58,0){$b_1$}
\psline[linewidth=.25](54,0)(54,1)(56,2)(54,3)(54,4)(52,5)(50,4)(50,0)
\end{pspicture}
\end{center}
We have the monomial pieces $t_k^{-1} \rightarrow {t'_i}^{-1} \&\ {t'_j}^{-1}t'_k$

%page 1 # 2

%page 2 # 2
\begin{center}
\begin{pspicture}(0,-.15)(20,1.9)
\psset{unit=.2}
\psline(0,10)(0,9)(2,8)(4,7)(4,6)(6,5)(8,4)(8,0)
\psline(4,10)(4,9)(2,8)(0,7)(0,3)(2,2)(4,1)(4,0)
\psline(8,10)(8,6)(6,5)(4,4)(4,3)(2,2)(0,1)(0,0)
\uput[90](0,10){$u_1$}
\uput[90](4,10){$u_2$}
\uput[90](8,10){$u_3$}
\uput[90](2,8){$t_i$}
\uput[90](6,5){$t_j$}
\uput[90](2,2){$t_k$}
\uput[270](0,0){$b_3$}
\uput[270](4,0){$b_2$}
\uput[270](8,0){$b_1$}
\psline[linewidth=.25](0,0)(0,1)(2,2)(0,3)(0,7)(2,8)(4,9)(4,10)

%1 212 piece

\psline(40,10)(40,6)(42,5)(44,4)(44,3)(46,2)(48,1)(48,0)
\psline(44,10)(44,9)(46,8)(48,7)(48,3)(46,2)(44,1)(44,0)
\psline(48,10)(48,9)(46,8)(44,7)(44,6)(42,5)(40,4)(40,0)
\uput[90](40,10){$u_1$}
\uput[90](44,10){$u_2$}
\uput[90](48,10){$u_3$}
\uput[90](46,8){$t_i'$}
\uput[90](42,5){$t_j'$}
\uput[90](46,2){$t_k'$}
\uput[270](40,0){$b_3$}
\uput[270](44,0){$b_2$}
\uput[270](48,0){$b_1$}
\psline[linewidth=.25](40,0)(40,4)(42,5)(44,6)(44,7)(46,8)(44,9)(44,10)

%2 212 piece
\psline(50,10)(50,6)(52,5)(54,4)(54,3)(56,2)(58,1)(58,0)
\psline(54,10)(54,9)(56,8)(58,7)(58,3)(56,2)(54,1)(54,0)
\psline(58,10)(58,9)(56,8)(54,7)(54,6)(52,5)(50,4)(50,0)
\uput[90](50,10){$u_1$}
\uput[90](54,10){$u_2$}
\uput[90](58,10){$u_3$}
\uput[90](56,8){$t_i'$}
\uput[90](52,5){$t_j'$}
\uput[90](56,2){$t_k'$}
\uput[270](50,0){$b_3$}
\uput[270](54,0){$b_2$}
\uput[270](58,0){$b_1$}
\psline[linewidth=.25](50,0)(50,4)(52,5)(54,4)(54,3)(56,2)(58,3)(58,7)(56,8)(54,9)(54,10)
\end{pspicture}
\end{center}
We have the monomial piece $t_k^{-1} \rightarrow {t'_i}^{-1} \&\ {t'_j}^{-1}t'_k$

% page 3 # 2
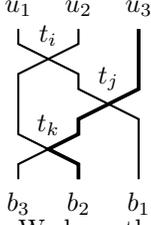
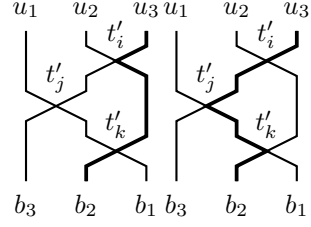
\begin{center}
\begin{pspicture}(0,-.15)(20,1.9)
\psset{unit=.2}
\psline(0,10)(0,9)(2,8)(4,7)(4,6)(6,5)(8,4)(8,0)
\psline(4,10)(4,9)(2,8)(0,7)(0,3)(2,2)(4,1)(4,0)
\psline(8,10)(8,6)(6,5)(4,4)(4,3)(2,2)(0,1)(0,0)
\uput[90](0,10){$u_1$}
\uput[90](4,10){$u_2$}
\uput[90](8,10){$u_3$}
\uput[90](2,8){$t_i$}
\uput[90](6,5){$t_j$}
\uput[90](2,2){$t_k$}
\uput[270](0,0){$b_3$}
\uput[270](4,0){$b_2$}
\uput[270](8,0){$b_1$}
\psline[linewidth=.25](4,0)(4,1)(2,2)(4,3)(4,4)(6,5)(8,6)(8,10)

%1 212 piece
\psline(40,10)(40,6)(42,5)(44,4)(44,3)(46,2)(48,1)(48,0)
\psline(44,10)(44,9)(46,8)(48,7)(48,3)(46,2)(44,1)(44,0)
\psline(48,10)(48,9)(46,8)(44,7)(44,6)(42,5)(40,4)(40,0)
\uput[90](40,10){$u_1$}
\uput[90](44,10){$u_2$}
\uput[90](48,10){$u_3$}
\uput[90](46,8){$t_i'$}
\uput[90](42,5){$t_j'$}
\uput[90](46,2){$t_k'$}
\uput[270](40,0){$b_3$}
\uput[270](44,0){$b_2$}
\uput[270](48,0){$b_1$}
\psline[linewidth=.25](44,0)(44,1)(48,3)(48,7)(46,8)(48,9)(48,10)

%2 212 piece
\psline(50,10)(50,6)(52,5)(54,4)(54,3)(56,2)(58,1)(58,0)
\psline(54,10)(54,9)(56,8)(58,7)(58,3)(56,2)(54,1)(54,0)
\psline(58,10)(58,9)(56,8)(54,7)(54,6)(52,5)(50,4)(50,0)
\uput[90](50,10){$u_1$}
\uput[90](54,10){$u_2$}
\uput[90](58,10){$u_3$}
\uput[90](56,8){$t_i'$}
\uput[90](52,5){$t_j'$}
\uput[90](56,2){$t_k'$}
\uput[270](50,0){$b_3$}
\uput[270](54,0){$b_2$}
\uput[270](58,0){$b_1$}
\psline[linewidth=.25](54,0)(54,1)(56,2)(54,3)(54,4)(52,5)(54,6)(54,7)(56,8)(58,9)(58,10)
\end{pspicture}
\end{center}
We have the monomial pieces  $t_k^{-1} \rightarrow {t'_i}^{-1} \&\ {t'_j}^{-1}t'_k$.

We can construct the other 16 cases of entry and exit points from these six using symmetry.  For each 
rigorous path on $R'$ listed above there is a corresponding rigorous path on $R$ as follows.  For the entry and exit 
points of the new path, exchange $b_1$ and $b_3$ (resp $u_1$ and $u_3$).  Then construct a new path which changes 
strings at the same vertices as the path on $R'$.  Unless the path on $R'$ changes strings at every vertex, this
new path will be rigorous.  This is what we call horizontal path symmetry. Rigorous paths which changes strings at 
every vertex do not have horizontally symmetric partners.  There is also a vertical symmetry.  Take a path on $R$ 
and exchange the ``top'' with ``bottom'' on the entry and exit points (i.e. $b_3$ is switched with $u_1$).  
Now create a path which changes strings at vertically flipped vertices (i.e. $t_i$ instead of $t_k$).
This will be a new rigorous path.  The remaining 16 single entry point cases can be gotten from the 6 listed above, 
through a combination of horizontal and vertical symmetries.  From this classification we see that the sets of
monomial pieces stated in the theorem are exchanged under a 3-move. 
\end{proof}

\end{subsection}
\begin{subsection}{Integrality of $\Delta^*$}
\begin{lemma}For any braid region $R$, the line connecting the apex of the $\lambda$-cone and the origin lies on the 
hyperplane $x_i+x_k=x_j$ .
\end{lemma}
\begin{proof}
It suffices to show that the apex of the $\lambda$-cone lies on the hyperplane.  At the apex all of the 
$\lambda$-inequalities are equalities. We can write the following
equations, where $A$, $B$ and are the contributions to the equation from $x_d$ for $d>k$.  
\begin{align*}x_i =& 2 + x_j - 2x_k + A \\
x_k =& 2 + B
\end{align*} Note that since $t_i$ and $t_k$ are in
the same column of the string diagram we have $A=B$.  From this we can see that 
\[x_i+x_k = (2+x_j-2x_k+A)+x_k=(x_j-x_k)+2+A=(x_j-x_k)+x_k=x_j.\]
\end{proof}
\begin{theorem}
$\Delta^*(\overline{\omega_0})$ is an integral polytope for any reduced decomposition $\overline{\omega_0}$ .
\end{theorem}
\begin{proof}
In \protect\cite{partialflag}, Batyrev (et. al.) show that for the reduced decomposition $\overline{\omega_0}$ = 
$s_1s_2s_1 \dots s_ns_{n-1} \dots s_1$ the polytope $\Delta^*$ is an integral polytope with a unique interior 
point $P=(1,2,1,\cdots,n,n-1,\cdots,2,1)$.  The facets of $\Delta$ determine linear functions $L$ such 
that $L(P)=1$.  Let $\Psi$ be the piecewise linear map between polytopes defined be Berenstein and Zelevinsky.  
We want to show that under a 3-move, the facets of $\Delta'$ determine $L'$ such that $L'(P')=1$ 
where $P' = \Psi(P)$.

The map $\Psi$ is defined by two linear maps $\psi_1$ and $\psi_2$ which agree on the hyperplane $x_i+x_k=x_j$.
For any linear function $L$ associated to a facet of $\Delta$ there are two corresponding linear functions on $\Delta'$;
$L_1'=L\circ \psi_1^{-1}$ and $L_2'=L\circ \psi_2^{-1}$.  

The point $P$ is halfway between the origin (which is the
apex of the string cone) and the apex of the $\lambda$-cone.  By the preceding lemma this shows that $P$ lies
on the hyperplane $x_i+x_k=x_j$.   Since $\psi_1$ and $\psi_2$ agree on that hyperplane we have $\psi_1(P)=\psi_2(P)=\Psi(P)$.

Therefore, we can make the calculations $L_1'(P')=L({\psi_1}^{-1}(P')=L(P)=1$ and $L_2'(P')=L({\psi_1}^{-1}(P')=L(P)=1$, which
implies that the polytope $\Delta'^{*}$ is integral.  Note that $P'$ is still halfway between the origin and the apex of the
$\lambda$-cone, which allows this argument to be repeated.

Under a 2-move the polytopes $\Delta$ and $\Delta'$ are isomorphic, so $\Delta^*$ is integral, if and only if $\Delta'^*$ is integral.  

By starting with $\overline{\omega_0}$ = $s_1s_2s_1 \dots s_ns_{n-1} \dots s_1$ and the point $P = 
(1,2,1,\cdots,n,n-1,\cdots,2,1)$ we can compose  braid moves and always have an interior point $P'$ on which all of the
defining linear functionals take value $1$.  This proves that $\Delta^*$ is integral for any choice of reduced decomposition.
\end{proof}
\begin{corollary}
The monomials in $f_{\overline{\omega_0},a}$ correspond to vertices of $\Delta^*$ and,
$F_{\overline{\omega_0}}$ is a subfamily of the linear system $|\OO_{X_{\Delta^*}}(1)|$ intersected with the torus.
\end{corollary}
\end{subsection}

\section{Birationality of $F_{\overline{\omega_0}}$}
\label{sec:birational}
\begin{definition}
Assume that $\overline{\omega_0}$ and $\overline{\omega_0}'$ differ by a 3-move.
Let $C = \frac{a_{\lambda_k}}{a_{\lambda_i}}$, $T$=Spec$\bC[t_1,t_1^{-1},\dots,t_N,t_N^{-1}]$
and $T'$=Spec$\bC[t'_1,{t'_1}^{-1},\dots,t'_N,{t'_N}^{-1}]$.  For a fixed coefficient vector $a$, define
$h_a : T' \dra T$ by \\
$(t_i,t_j,t_k) = (\frac{t_i't_k' + Ct_j'}{t_i'},\frac{t_i't_k'}{C},\frac{Ct_i't_j'}{t_i't_k'+Ct_j'}) \text{, and }
t_q = t_q' \text{ for }q \notin \{i,j,k\}$.  \\

\end{definition}
Define a map $h_a': T \dra T'$ by exchanging $t$ with $t'$ in the map above.  Then we have $h_a' \circ h_a = id.$
\begin{proposition}
\label{prop:phi}
For a fixed coefficient vector $a=(a_1,\cdots,a_r) \in P_{\overline{\omega_0}}$,  $f' := h_a^*(f_{\overline{\omega_0},a})$ defines a 
variety $Z_{\overline{\omega_0}',a'}$ for some coefficient vector $a'$.
Therefore, we have an induced birational map $h_a: Z_{\overline{\omega_0}',a'} \dra Z_{\overline{\omega_0},a}$.
\end{proposition}
\begin{remark}
\label{rem:g}
We will see in Proposition \protect\ref{prop:iso} that the map $g$ defined by $g^*(a)=a'$ is an isomorphism between
the parameter spaces $P'_{\overline{\omega_0}'}$ and $P_{\overline{\omega_0}}$.  The coordinates for $a'$ are explicitly given as regular functions of the
coordinates for $a$ in Proposition \protect\ref{prop:monomials}.
\end{remark}
\begin{remark}
The map $h_a$ is a specific geometric lift (in the sense of \protect\cite{bz}) of the piecewise linear function: 
($t_i',t_j',t_k'$) $\rightarrow$ (min($t_k'$,$t_j' - t_i'$),$t_i' + t_k'$,max($t_i'$,$t_j' - t_k'$)). 
\end{remark}
\emph{Proof of Proposition \protect\ref{prop:phi}.} 
We show that for a fixed coefficient vector $a$, $h_a$ is birational in two steps.  First, in Proposition \ref{prop:monomials}, we show 
that $h_a^*(1 - \sum_{d} a_iT^{M_d})= (1 - \sum_{d'} a_i'T^{M_d'})$ for some coefficient vector $a'$.
Then, in Section \ref{sec:global} we show that the coefficient vector $a'$ is actually in $P'_{\overline{\omega_0}'}$ by showing that $g^*$ preserves the combinatorial box equations.  
\begin{subsection}{Action of $h^*$ on $f_{\overline{\omega_0},a}$}
\label{sec:local}
\begin{definition}
We refer to the box equations coming from the boxes $R$ and $R'$ as \emph{local box equations}.
\end{definition}
\begin{proposition}
\label{prop:monomials}
Fix a coefficient vector $a$ satisfying the local box equations.  Then 
$h_a^*(1 - \sum_{d} a_dT^{M_d})= (1 - \sum_{d'} a_d'T^{M_d'})$ for some coefficient vector $a'$.
\end{proposition}
\begin{proof}
We examine how $h_a^*$ acts on $1 - \sum_{d} a_iT^{M_d}$ monomial by monomial. Since $h_a^*$ only depends on $t_i,t_j,$ and $t_k$
we classify its action on monomials based on the corresponding monomial pieces.  

A priori  $h_a^*(f_{\overline{\omega_0},a})$ is  a rational function, but not a sum of monomials. 
In several cases, in order to see that $h_a^*(1 - \sum_{d} a_iT^{M_d})$ is a linear combination of monomials, we have to group pairs of monomials together.  
In some cases $h_a^*$ only maps a particular combination
of two monomials to a monomial, if the coefficients satisfy the local box equations.  \\
The map $h_a^*$ take the following classes of monomials to monomials (or sum of monomials):
\begin{align*}
 h_a^*(a_dt_i^{-1}t_k^{-1}) &= a_da_{\lambda_i}a_{\lambda_k}^{-1}t_j'^{-1} \\
 h_a^*(a_dt_i^2t_j^{-1}t_k^2) &= a_da_{\lambda_i}^{-3}a_{\lambda_k}^3t_i'^{-1}t_j'^2t_k'^{-1}  \\
 h_a^*(a_dt_i^{-1}t_j^2t_k^{-1}) &= a_da_{\lambda_i}^3a_{\lambda_k}^{-3}t_i'^2t_j'^{-1}t_k'^2  \\
 h_a^*(a_dt_j^{-1}) &= a_da_{\lambda_i}^{-1}a_{\lambda_k} t_i'^{-1}t_k'^{-1} \\
 h_a^*(a_dt_i^{-1}t_k^{-1}) &= a_da_{\lambda_i}a_{\lambda_k}^{-1}t_j'^{-1}  \\
 h_a^*(a_dt_j^{-1}) &= a_da_{\lambda_i}^{-1}a_{\lambda_k}t_i'^{-1}t_k'^{-1}  \\
 h_a^*(a_dt_it_k) &= a_da_{\lambda_i}^{-1}a_{\lambda_k}t_j'  \\
 h_a^*(a_dt_j) &= a_da_{\lambda_i}a_{\lambda_k}^{-1}t_i't_k'   \\
 h_a^*(a_dt_it_j^{-1}t_k) &= a_da_{\lambda_k}^2a_{\lambda_i}^{-2}t_i'^{-1}t_j't_k'^{-1} \\
 h_a^*(a_{\lambda_j}t_jt_k^{-1}) &= a_{\lambda_i}^2a_{\lambda_j}a_{\lambda_k}^{-2}t_i't_j'^{-1}t_k'^2 + a_{\lambda_i}a_{\lambda_j}a_{\lambda_k}^{-1}t_k'  \\
 h_a^*(a_dt_i^{-1}t_jt_k^{-1}) &= a_da_{\lambda_k}^{-2}a_{\lambda_i}^2t_i't_j'^{-1}t_k' \\
 h_a^*(a_dt_k^{-1}) &= a_dt_i'^{-1} + a_da_{\lambda_i}a_{\lambda_k}^{-1}t_j'^{-1}t_k'  \\
 h_a^*(a_dt_i) &= a_dt_k' + a_da_{\lambda_i}^{-1}a_{\lambda_k}t_i'^{-1}t_j'  \\
\end{align*}
The following classes of monomials must be grouped together in order for $h_a^*$ to take them to a monomial.
\begin{align*}
 h_a^*(a_{\lambda_i}t_it_j^{-1}t_k^2 + a_{\lambda_k}t_k) &= a_{\lambda_i}^{-1}a_{\lambda_k}^2t_j't_k'^{-1}  \\
\end{align*}
The following classes of monomials must be grouped together, and their coefficients must satisfy the local box equations
in order for $h_a^*$ to take them to monomials.  We write the
second coefficient in terms of the local box equation.  
\begin{align*}
 h_a^*(a_dt_k+a_da_{\lambda_i}^{-1}a_{\lambda_k}t_i^{-1}t_j) &= a_dt_i'  \\
 h_a^*(a_dt_i^{-1} + a_da_{\lambda_i}a_{\lambda_k}^{-1}t_j^{-1}t_k) &= a_dt_k'^{-1}  
\end{align*}
The monomials occurring in $\image({h_a}^*)$ are precisely those occurring in $f_{\overline{\omega_0},a}$.
From this we see that \\ $h_a^*(1 - \sum_{d} a_iT^{M_d})= (1 - \sum_{d'} a_i'T^{M_d'})$ for some coefficient vector $a'_i$.  
The inverse map is given by using the same construction for the braid move
from $\overline{\omega_0}'$ to $\overline{\omega_0}$.
\end{proof}
\end{subsection}
\begin{subsection}{Preservation of Box Equations.}
\label{sec:global}
\begin{proposition}
For any 3-move, the map $g^*$ preserves box equations.
\end{proposition}
\begin{proof}
Let $R'$ be the image of the box $R$  (We can assume it is of type $R212$ with the proof of the other
case following using the same methods). We classify all boxes $O'$ of the string diagram for
$\overline{\omega_0}'$ by their positions in comparison with $R'$. 

The proof of every case follows the following method:  First, classify all possibilities of monomial
pieces $\tau$ for which the box condition on $O'$ could be satisfied.  Since we know how
the monomial pieces change under the braid move, we verify that the box conditions must be
satisfied for the box $O=\preimage(O')$.  Since $a \in P$ we know that if the box conditions are satisfied on $O$, then
the corresponding box equations (in $a_d$) are also satisfied. We use the map $g$ to write $a_d'$ in terms of $a_d$,
and check that  new combinatorial box equations on $O'$ are satisfied. 

\begin{remark}
For some string monomials on $O'$ there are two different corresponding
string monomials in the preimage.  Only one of these two string monomials is used to construct a
box condition on the preimage. The first box 
condition may not be met in the preimage when the other choice of string monomial is used.
\end{remark}
We classify the boxes $O'$ into eight groups.  For each group, we state how many pair of monomial pieces could
possibly correspond to monomials satisfying the box conditions on $O'$.
\begin{enumerate}
\item $O'$ doesn't touch $R'$ (1 pair)
\item $O' = R'$ (1 pair)
\item $O'$ is below and left of $R'$ (3 pair)
\item $O'$ is above and left of $R'$ (3 pair)
\item $O'$ is directly above $R'$ (3 pair)
\item $O'$ is directly below $R'$ (3 pair)
\item $O'$ is to the right of $R'$ (6 pair)
\item $O'$ is to the left of $R$ (6 pair)
\end{enumerate}
What follows is the proof for the case when $R'$ is directly above $O'$.  The proof for all of the other groups follow the exact
same method.\\

Assume $R'$ is directly above $O'$ as pictured below.
% case 3
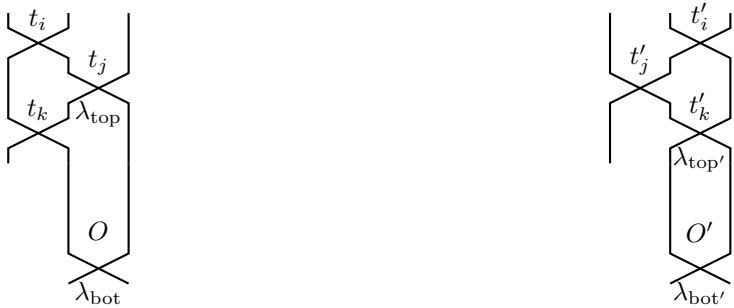
\begin{center}
\begin{pspicture}(0,-2)(20,2.2)
\psset{unit=.2}
\psline(0,10)(0,9)(2,8)(4,7)(4,6)(6,5)(8,4)(8,0)
\psline(4,10)(4,9)(2,8)(0,7)(0,3)(2,2)(4,1)(4,0)
\psline(8,10)(8,6)(6,5)(4,4)(4,3)(2,2)(0,1)(0,0)

\uput[90](2,8){$t_i$}
\uput[90](6,5){$t_j$}
\uput[90](2,2){$t_k$}
\uput[270](6,5){$\lambda_{\TOP}$}
\uput[270](6,-7){$\lambda_{\bottom}$}
\uput[270](6,-3){$O$}
\psline(8,0)(8,-6)(6,-7)(4,-8)
\psline(4,0)(4,-6)(6,-7)(8,-8)

%1 212 piece
\psline(40,10)(40,6)(42,5)(44,4)(44,3)(46,2)(48,1)(48,0)
\psline(44,10)(44,9)(46,8)(48,7)(48,3)(46,2)(44,1)(44,0)
\psline(48,10)(48,9)(46,8)(44,7)(44,6)(42,5)(40,4)(40,0)

\uput[90](46,8){$t_i'$}
\uput[90](42,5){$t_j'$}
\uput[90](46,2){$t_k'$}
\uput[270](46,2){$\lambda_{\TOP'}$}
\uput[270](46,-7){$\lambda_{\bottom'}$}
\uput[270](46,-3){$O'$}
\psline(48,0)(48,-6)(46,-7)(44,-8)
\psline(44,0)(44,-6)(46,-7)(48,-8)
\end{pspicture}
\end{center}

Assume that $p_1'$ and $p_2'$ satisfy the box conditions on $O'$.  Let $\tau^{M_d}$ be the monomial piece associated to
$T^{M_d}$.  By comparing the $\lambda$-inequalities see the following relationship among monomial pieces.  
\[\frac{\tau^{p_2'}}{\tau^{p_1'}} = \frac{\tau^{\lambda_{\TOP'}}}{\tau^{\lambda_{\bottom'}}} = t_k.\] \\
From this we see that there are three possibilities for pairs string monomial pieces.
\begin{enumerate}
\item Assume that
\[
\tau^{p_1'} = t_i^{-1}t_k^{-1} \text{ and }
\tau^{p_2'}=t_i^{-1}.
\]
We construct a set of monomials  satisfying the box conditions for
the box $O$ as follows. 
\begin{align*}
T^{p_1'}  \cdot  T^{\lambda_{\TOP'}}  &=  T^{p_2'}  \cdot  T^{\lambda_{\bottom'}}  \\
T^{p_1'}  \cdot  T^{\lambda_{\TOP}} \cdot t_j^{-1}t_k^2  &=  T^{p_2'}  \cdot  T^{\lambda_{\bottom}}  \\
T^{p_1} \cdot t_i^{-1}t_jt_k^{-1}  \cdot  T^{\lambda_{\TOP}} \cdot t_j^{-1}t_k^2  &=  T^{p_2} \cdot t_i^{-1}t_k \cdot  T^{\lambda_{\bottom}}  \\
T^{p_1} \cdot T^{\lambda_{\TOP}}  &=  T^{p_2} \cdot  T^{\lambda_{\bottom}}  
\end{align*}
Since the box condition are satisfied on $O$, we have the following  equations in $a_i.$ \\
\[
a_{p_1}a_{\lambda_{\TOP}} = a_{p_2}a_{\lambda_{\bottom}}.\]   

Rewriting the $a_i'$  in terms of $a_i$ we see that \\
\[
a_{p_1'}a_{\lambda_{\TOP'}} = a_{p_1} a_{\lambda_{\TOP}} \text{ and } 
a_{p_2'}a_{\lambda_{\bottom'}} = a_{p_2}a_{\lambda_{\bottom}} . \]\\
Therefore, 
\[a_{p_1'}a_{\lambda_{\TOP'}} = a_{p_2'}a_{\lambda_{\bottom'}}. \]
\item Assume that
\[\tau^{p_1'} = t_k^{-1} \text{ and } \tau^{p_2'}=1.\] \\
We construct a set of monomials  satisfying the box conditions for
the box $O$ as follows. 
\begin{align*}
 T^{p_1'}  \cdot  T^{\lambda_{\TOP'}}  &=  T^{p_2'}  \cdot  T^{\lambda_{\bottom'}}  \\
 T^{p_1'}  \cdot  T^{\lambda_{\TOP}} \cdot t_j^{-1}t_k^2  &=  T^{p_2'}  \cdot  T^{\lambda_{\bottom}}  \\
 T^{p_1}  \cdot  t_jt_k^{-2}  \cdot  T^{\lambda_{\TOP}} \cdot t_j^{-1}t_k^2  &=  T^{p_2}  \cdot  T^{\lambda_{\bottom}} \\
 T^{p_1} \cdot T^{\lambda_{\TOP}}  &=  T^{p_2} \cdot  T^{\lambda_{\bottom}} 
\end{align*} 
Since the box condition were satisfied on $O$, we have the following equations in $a_i.$ \\
\[ 
 a_{p_1}a_{\lambda_{\TOP}} = a_{p_2}a_{\lambda_{\bottom}} .\] \\
Rewriting the $a_i'$  in terms of $ a_i$  we see that 
\[a_{p_1'}a_{\lambda_{\TOP'}} = a_{p_1}a_{\lambda_{\TOP}} 
 \text{ and }
a_{p_2'}a_{\lambda_{\bottom'}} = a_{p_2}a_{\lambda_{\bottom}}. \]
Therefore, 
\[a_{p_1'}a_{\lambda_{\TOP'}} = a_{p_2'}a_{\lambda_{\bottom'}}.\]
\item Assume that
\[\tau^{p_1'} = t_i^{-1}t_jt_k^{-1} \text{ and } \tau^{p_2'}=t_i^{-1}t_j.\] \\
We construct a set of monomials satisfying the box conditions for
the box $O$ as follows. 
\begin{align*}
 T^{p_1'}  \cdot  T^{\lambda_{\TOP'}}  &=  T^{p_2'}  \cdot  T^{\lambda_{\bottom'}}  \\
 T^{p_1'}  \cdot  T^{\lambda_{\TOP}} \cdot t_j^{-1}t_k^2  &=  T^{p_2'}  \cdot  T^{\lambda_{\bottom}}  \\
 T^{p_1}  \cdot  t_i^{-2}t_j \cdot  T^{\lambda_{\TOP}} \cdot t_j^{-1}t_k^2  &=  T^{p_2} t_i^{2}t_j^{-2}t_k^{2} \cdot  T^{\lambda_{\bottom}} \\
 T^{p_1} \cdot T^{\lambda_{\TOP}}  &=  T^{p_2} \cdot  T^{\lambda_{\bottom}} 
\end{align*} 
Since the box condition are satisfied on $O$, we have the following equations in $a_i.$ \\
\[
 a_{p_1}a_{\lambda_{\TOP}} = a_{p_2}a_{\lambda_{\bottom}} .\] \\
Rewriting the $a_i'$  in terms of $ a_i$  we see that 
\[a_{p_1'}a_{\lambda_{\TOP'}} = a_{p_1}\frac{a_k^2}{a_i^2}a_{\lambda_{\TOP}}\frac{a_i}{a_k} 
\text{ and } 
a_{p_2'}a_{\lambda_{\bottom'}} = a_{p_2}\frac{a_k}{a_i}a_{\lambda_{\bottom}}.\]
Therefore, 
\[a_{p_1'}a_{\lambda_{\TOP'}} = a_{p_2'}a_{\lambda_{\bottom'}}.\]
\end{enumerate}
We have shown that the combinatorial box equations are preserved under $g^*$, which completes the proof of
proposition \protect\ref{prop:phi}.
\end{proof}
\end{subsection}
\begin{subsection}{Isomorphism of Parameter Spaces}
\begin{proposition}
\label{prop:iso}
If $\overline{\omega_0}$ and $\overline{\omega_0}'$ differ by a single braid move then $P'_{\overline{\omega_0}'} \cong P_{\overline{\omega_0}}$.
\end{proposition}
\begin{proof}
Assume $\overline{\omega_0}$ and $\overline{\omega_0}'$ differ by a 3-move.  As in remark \protect\ref{rem:g}, define $g:P' \rightarrow P$ by $g^*(a)=a'$.  
By the proof of Proposition 4.6 we see that this defines a regular map from $P_{\overline{\omega_0}}$ to the
torus of coordinates of $P'_{\overline{\omega_0}'}$.  Proposition 4.7 shows that the image of this map is actually in $P'_{\overline{\omega_0}'}$.
Similarly we can construct a regular map $g': P_{\overline{\omega_0}} \rightarrow P'_{\overline{\omega_0}'}$.
We can use Proposition \protect\ref{prop:monomials} to verify that these maps are inverse to one another.  

To complete the proof we note that $P_{\overline{\omega_0}}$ and $P'_{\overline{\omega_0}'}$ are obviously isomorphic under a 2-move. 
\end{proof}
\begin{corollary}
\label{lemma:irred}
The parameter space $P_{\overline{\omega_0}}$ is a toric variety for any choice of $\overline{\omega_0}$.
\end{corollary}
\begin{proof}
In the case of the standard reduced decomposition ($s_1s_2s_1 \dots s_ns_{n-1} \dots s_1$) the combinatorial box equations
and Batyrev's box equations are exactly the same.  By Batyrev's work \protect\cite[Rm.4.2]{B_mirror}, the set of non-zero $a_i$ satisfying
the box equations form a toric variety.  By the previous lemma, we see that under a braid move
the parameter spaces are isomorphic.  We can compose braid moves to see that the parameter spaces are isomorphic for any choice
of reduced decomposition.  
\end{proof}
\end{subsection}
\begin{subsection}{Birationality of the Families $F_{\overline{\omega_0}}$}
\begin{lemma}
For $\overline{\omega_0}$ and $\overline{\omega_0}'$ differing by a 2-move, the map $H: F_{\overline{\omega_0}'} \rightarrow F_{\overline{\omega_0}}$
defined by exchanging $t_i$ with $t_j$ and exchanging the corresponding coefficients is an isomorphism.
\end{lemma}
\begin{proof}
This is a consequence of Berenstein-Zelevinsky's map between string cones for the 2-move case.
\end{proof}
\begin{theorem}
If $\overline{\omega_0}$ and $\overline{\omega_0}'$ differ by
a 3-move, the map \[H: F_{\overline{\omega_0}'} \dra F_{\overline{\omega_0}}\] defined as $H:=(g,h_a)$
is a birational map of families.
\end{theorem}
\begin{proof}
This follows from the fact that $g$ is an isomorphism, and $h_a$ is a birational map for any fixed coefficient vector $a$.
\end{proof}
\begin{corollary}
\label{cor:birational}
For any two reduced decompositions $\overline{\omega_0}$ and $\overline{\omega_0}'$, the families $F_{\overline{\omega_0}}$ and 
$F_{\overline{\omega_0}'}$ are birational.
\end{corollary}
\begin{proof}
Connect $\overline{\omega_0}$ and $\overline{\omega_0}'$ by a sequence of braid moves.  For each braid move
there is a birational map of families, so we can compose these maps to get a birational map between $F_{\overline{\omega_0}}$ and 
$F_{\overline{\omega_0}'}$. 
\end{proof}
\begin{corollary}
Smooth Calabi-Yau manifolds birational to the generic elements of $F_{\overline{\omega_0}}$ have the
same Hodge numbers for any choice of $\overline{\omega_0}$.
\end{corollary}
\begin{proof}
Since the generic elements are all birational, smooth Calabi-Yau manifolds birational to them have the same Hodge numbers 
by the work of Batyrev and Kontsevich \protect\cite{bbetti, Kont}.
\end{proof}
\end{subsection}
\section{Comparing $F_{\overline{\omega_0}}$ with the Families Defined by Batyrev (et.al.)}
\label{sec:small}
\begin{definition}\protect\cite[Def 3.1]{B_mirror}
\label{Def:small}
Let $X \subset$ $\bP^m$ be a smooth Fano variety of dimension $n$.  A normal Gorenstein toric Fano variety
$Y \subset \bP^m$ is called a \emph{small toric degeneration} of $X$, if there exists a Zariski open neighbourhood
$U$ of $0 \subset \bA^1$ and an irreducible subvariety $\wX \subset \bP^m \times U$ such that the morphism
$\pi : \wX \mapsto U$ is flat and the following conditions hold:
\begin{enumerate}
\item the fiber $X_t := \pi^{-1}(t) \subset \bP^m$ is smooth for all $t \in U \backslash 0$;
\item  the special fiber $X_0 := \pi^{-1}(0) \subset \bP^m$ has at worst Gorenstein terminal singularities
and $X_0$ is isomorphic to $Y \subset \bP^m$;
\item the canonical homomorphism \[\Pic(\wX /U) \mapsto \Pic(X_t) \] is an isomorphism for all $t \in U$.
\end{enumerate}
\end{definition}
In \cite{B_mirror} Batyrev proposes a mirror construction for any small toric degeneration. This construction is 
an extension of Batyrev (et. al.)'s construction in the case of partial flag varieties.  We set up the next lemma to help
show that our construction is the same as Batyrev's in the case of a small degeneration.
\begin{lemma}
\label{lemma}
For any inequalities $\lambda_{\TOP},\lambda_{\bottom},p_1$, and $p_2$ whose corresponding monomials satisfy a combinatorial box condition on O,
there exists a facet of $\Delta^*$ 
containing the vertices corresponding to $\lambda_{\TOP},\lambda_{\bottom},p_1$, and $p_2$.
\end{lemma}
\begin{proof}
Define the linear function $L(t_1, \cdots,t_N) := t_{\TOP} - t_{\bottom} + 1$.
We check that $L$ = $0$ on the corresponding vertices $\lambda_{\TOP},\lambda_{\bottom},p_1$, and $p_2$, and $L \ge 0$
on all other vertices of $\Delta^*$.\\
\begin{align*}
L(\lambda_{\TOP}) &= 1 - 2 + 1 = 0 \\
L(\lambda_{\bottom}) &= 0 -1 + 1 = 0 \\
L(p_1) &= -1 - x + 1 \\
L(p_2) &= y - 1 + 1 \\
\end{align*} 
But since $\lambda_{\TOP} + p_1 = \lambda_{\bottom} + p_2$ we see that $x=y=0$ which implies that
$L(p_1)=L(p_2)=0$.

Next we want to show that for every other vertex $v$ of $\Delta^*, L(v) \ge 0$.
Assume $v$ is a vertex corresponding to a $\lambda$-inequality $\lambda_v$ then \\
\[L(v)=
\begin{cases}
0  - 0    + 1 &= 1 \text{ If $\lambda_v$ lies below $\lambda_{\bottom}$ or least 2 columns to the} \\ 
&  \text{left or right of $\lambda_{\TOP}$} \\

0  -(-1)  + 1 &= 2 \text { If $\lambda_v$ is in the region O.} \\
2  - 2    + 1 &= 1 \text{ If $\lambda_v$ lies in the column above $\lambda_{\TOP}$} \\
-1 - (-1) + 1 &= 1 \text{ If $\lambda_v$ lies above and adjacent to  $\lambda_{\TOP}$.} \\
\end{cases}
\]
String vertices  can only take values from the set \{-1,0,1\}
on the coordinates $t_{\TOP}$ and $t_{\bottom}$.  From this we see that the only way the linear function
$L$ could be negative on such a vertex $v$, is if  $t_{\TOP}$ coordinate of $v$ was $-1$ and the $t_{\bottom}$ coordinate was $1$.  
Any such inequality would correspond to an oriented path which intersects itself, which is a contradiction.

We see that $(L=0)$ defines a face of $\Delta$ containing the four vertices.  This completes the proof since any face
is contained in a facet.
\end{proof}

\begin{theorem}
\label{thm:small} For a small string degeneration $F_{\overline{\omega_0}}$ 
is  same as the family constructed by Batyrev in \protect\cite[Sect. 4]{B_mirror}.
\end{theorem}
\begin{proof}
In the case of a small toric degeneration Batyrev defines a subfamily of varieties in $|\OO_{X_{\Delta^*}}(1)|$
whose coefficients satisfy box equations. (Box equation terminology was used in \protect\cite{partialflag}.  The coefficients
satisfying the box equations are referred to as $\Sigma$-admissible in \protect\cite{B_mirror}.)

Our combinatorial box equations are a subset of Batyrev's
box equations if, and only if, the vertices in $\Delta^*$ corresponding to $\lambda_{\TOP},\lambda_{\bottom},p_1$, and $p_2$
all lie in a common facet of $\Delta^*$.
Therefore, by Lemma \protect\ref{lemma} we have shown that the combinatorial box equations are a subset of Batyrev's box equations.  
Batyrev proves in \protect\cite[Rem. 4.2]{B_mirror}  that his box equations define an irreducible variety
of dim $= \rank(\Pic(X_{\Delta^*})) + \dim(G/B)$.  This is the same as the
dimension of the parameter space $P_{\overline{\omega_0}}$.  By Lemma \protect\ref{lemma:irred}, $P_{\overline{\omega_0}}$ is an irreducible variety. 
If Batyrev's family had any more independent box equations they would define a variety of smaller dimension. 
Therefore, Batyrev's box equations cut out exactly $P_{\overline{\omega_0}}$, and the families are the same.
\end{proof}
\begin{corollary}
For all reduced decompositions $\overline{\omega_0}$, generic elements of $F_{\overline{\omega_0}}$ share the mirror properties of the generic
elements of Batyrev's families.
\end{corollary}
\begin{proof}
In the small case the families are the same by Theorem \protect\ref{thm:small}, and thus their generic elements share the same mirror properties.  By
Corollary 4.12 the generic elements are birational for any choice of reduced decomposition, and therefore share
the same mirror properties.
\end{proof}
\begin{remark}
For the reduced decomposition $\overline{\omega_0} = s_3s_1s_2s_1s_3s_2$ calculation shows that rank($\Pic(X_{\Delta}$)) $<$ rank($\Pic(G/B$)).
Therefore, the corresponding string degenerations is not small. This calculation was performed with the
aid of the freely available program $\tt polymake$ using the techniques found in \protect\cite{fulton}.
\end{remark}

\bibliographystyle{amsalpha}
%\bibliography{primary}

\def\cprime{$'$}
\providecommand{\bysame}{\leavevmode\hbox to3em{\hrulefill}\thinspace}
\providecommand{\MR}{\relax\ifhmode\unskip\space\fi MR }
% \MRhref is called by the amsart/book/proc definition of \MR.
\providecommand{\MRhref}[2]{%
  \href{http://www.ams.org/mathscinet-getitem?mr=#1}{#2}
}
\providecommand{\href}[2]{#2}

\end{document}